%% file: main.tex
\documentclass[]{article}
\usepackage[pagebackref=false,
            colorlinks=true,
            linkcolor=black,
            bookmarks=true,
            filecolor=black,
            urlcolor=black,
            citecolor=black
           ]{hyperref}
\usepackage{arxiv}

\usepackage{graphicx}      
\usepackage{amsmath}
\usepackage{amssymb}
\usepackage{siunitx}
\usepackage[utf8]{inputenc} 
\usepackage[T1]{fontenc}    
\usepackage{booktabs}       
\usepackage{amsfonts}       
\usepackage{nicefrac}       
\usepackage{microtype}      
\usepackage{doi}
\usepackage{natbib}
\DeclareMathOperator*{\E}{\mathbb{E}}

\title{An Open Source Stochastic Unit Commitment Tool using the PyPSA-Framework}

\date{May 10, 2024}	

\author{Tom Welfonder \\
        University of Stuttgart \\
        Stuttgart, Germany \\
        \And
        Johannes Lips \\
	Institute of Combustion and Power Plant Technology\\
	University of Stuttgart\\
	Stuttgart, Germany \\
	\texttt{johannes.lips@ifk.uni-stuttgart.de} \\
        \And
        Alois Gmur \\
        University of Stuttgart \\
        Stuttgart, Germany \\
        \And
	\href{https://orcid.org/0000-0002-0208-4100}{\includegraphics[scale=0.06]{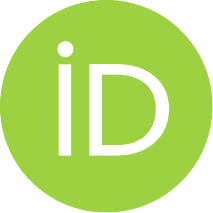}\hspace{1mm}Hendrik Lens} \\
	Institute of Combustion and Power Plant Technology\\
	University of Stuttgart\\
	Stuttgart, Germany \\
}

\hypersetup{
pdftitle={An Open Source Stochastic Unit Commitment Tool using the PyPSA-Framework},
pdfauthor={T. Welfonder, J. Lips, A. Gmur, H. Lens},
pdfkeywords={unit commitment, stochastic optimization, electricity markets, uncertainty modeling, open source, waste-to-energy, CHP, OR in energy},
}

\begin{document}
\maketitle

\begin{abstract}
\input{abstract}
\end{abstract}



\setcounter{footnote}{0} 
\section{Introduction}
For a portfolio of power generators and storage units that is managed collectively, unit commitment (UC) strategies, also called optimal dispatch strategies, may be used in order to maximize profits.
Through UC optimization, the portfolio manager tries to find the operating schedule of each portfolio component that will lead to maximal profit, taking into account various economical inputs and technical constraints.
The problem might, e.g., consider the (expected) electricity market prices, fuel and maintenance costs, as well as process efficiencies, minimal loads, and start-up and shut-down times.
Any product delivery agreement concerning electricity, industrial steam, or district heating or cooling can be considered through additional constraints for the optimization problem.
Mathematically, the UC problem is commonly formulated as a Mixed-Integer Linear Program (MILP).

Parallel to the increase of volatile renewable generation and battery energy storage systems (BESS), the problem has been expanded to include these technologies, and has been applied to virtual power plants.
To deal with the uncertainty of (volatile) production capacity, demand, and price prognosis, stochastic UC optimization has been developed \citep{zheng_stochastic_2015}.

UC optimization is normally implemented using commercial or in-house proprietary tools.
To the best of our knowledge, no comparable open source tool is available that offers stochastic UC for a modular and easy-to-configure portfolio setup.
In this paper, we present an open source stochastic UC tool, which is available on GitHub\footnote{https://github.com/PPGS-Tools/PyPSA-stochUC} \citep{welfonder_pypsa-stochuc_2024} and present an example use case in which stochastic UC optimization is done for a waste-to-energy plant with heat storage and BESS in Germany.
Using the tool, it is possible to optimize the UC of a customizable portfolio for sales on the day-ahead (DA) market and aFRR (balancing power) markets, resulting in an hourly operating schedule for all portfolio components.
We hope the tool lowers the threshold for using UC optimization, especially for small portfolios, ranging from regional utility providers to families who own, e.g., electric vehicles, solar panels and a heat pump.
This would result in a more efficient use of resources.

The UC tool builds on the Python library PyPSA \citep{brown_pypsa_2018}.
PyPSA is a toolbox for energy system optimization and simulation that has grown significantly in functionality in recent years.
With PyPSA, it is possible to construct a network of high-level system components, translate it to an optimization problem and solve the problem.
PyPSA contains readily-available system components such as committable production units (`generators') for which nominal load, marginal costs, efficiency, minimal load and other specifications can be given.
Other components with similar complexity include `storage units', `loads', `buses' (to which all other components are connected), and `links' (which connect buses directly to each other with a given efficiency).
The strength of PyPSA is that its modular structure makes it easy to add and change variables and constraints, making it the ideal basis for the UC tool.

To make stochastic UC optimization possible in PyPSA, we have implemented multiple extensions.
This is done in a modular way, so that the extensions integrate naturally in PyPSA and can be used separately or combined.
In the discussion of the methodology, each extension is presented separately, making it easy for other developers to use that specific extension in other PyPSA projects according to their needs.
In the example that we present at the end of the paper, all extensions are used together in the stochastic UC optimization of a waste-to-energy plant with storage facilities.
The example shows that the UC optimization tool can be used successfully under the consideration of electricity prices and district heating demand uncertainty.

\section{Methodology}\label{sec:methods}
PyPSA is an energy system optimization (ESO) toolbox, normally used to model large power systems, such as the European power system, and perform optimal power flow or total least-cost optimizations on them.
PyPSA does not have built-in functionality for stochastic optimization, nor for market mechanisms.

The PyPSA network structure, model translation and least-cost optimization are used as the basis for the UC optimization.
Specifically, the portfolio that is to be optimized is constructed as a PyPSA network.
Several functionalities are added to PyPSA to enable stochastic UC optimization of this network:
\begin{enumerate}
    \item \textbf{Market Mechanisms} -- In UC optimization, only a small number of market participants are modelled. Under a price-taker assumption, market prices are inputs of the optimization. Market and bidding mechanisms should be part of the model.
    \item \textbf{Stochastic Optimization} -- Scenario-based stochastic optimization is used to obtain robust UC under uncertainty. In this context, a classical PyPSA optimization can be seen as single-scenario optimization. Implementing multiple-scenario optimization can be done without changing the basic framework of least-cost optimization in PyPSA.
    \item \textbf{Rolling Horizon and Multistaging} -- In UC, only partial information is available at each instance in time. Rolling horizon optimization is a common approach for dealing with forecasts that are more accurate the closer they come to realization. If multiple decisions, such as market bid submissions, need to be made at different times, a multistage optimization can be used in combination with a rolling horizon.
\end{enumerate}

The optimization problem is formulated analogously to \citet{kraft_stochastic_2022} and \citet{zheng_stochastic_2015} as a MILP problem.
The least-cost optimization of the portfolio network can be written in the following form, with $n_\mathrm{r}$ being the number of real decision variables, $n_\mathrm{b}$ the number of binary decision variables, $n=n_\mathrm{r}+n_\mathrm{b}$, and $m$ the number of constraints:
\begin{equation}\label{eq:MILP}
    \begin{aligned}
        \min_x{J(x)}&= c^\intercal x,  &c \in \mathbb{R}^{n \times 1}, x \in \mathbb{R}^{n_\mathrm{r} \times 1} \times \{0,1\}^{n_\mathrm{b}\times 1}, \\
        \text{s.t. } A x &\leq b, &A \in \mathbb{R}^{m \times n}, b \in \mathbb{R}^{m \times 1}.
    \end{aligned}
\end{equation}
Here, $x$ are the decision variables, which contain the operating schedules of all portfolio components, as well as energy storage levels.
$c$ are the costs related to $x$ in the objective function $J(x)$, corresponding with the marginal costs associated with fuel or maintenance.
$A \, x \leq b$ formulate constraints that apply to the problem.

\subsection{Power and Energy Markets}\label{sec:markets}
We consider two markets: the Day-Ahead Electricity Market (DAM) and the Automatic Frequency Restoration Reserve (aFRR) Market.
aFRR is one of the balancing energy products in Europe for which there is a market, mostly in the form of an auction, in each country.
Other markets, continuous intraday trading, and aFRR activation are not yet considered in the context of this paper, as they would render the optimization problem significantly more complex.
For both markets, the physical implications of trading as well as the bidding and market mechanisms are modelled according to their implementation in Germany.
The model formulations are based on \citet{kraft_stochastic_2022} and \citet{kumbartzky_optimal_2017}, which are adapted to be compatible with PyPSA.

\subsubsection{Day-Ahead Market}
A PyPSA-`generator', representing the grid to which the portfolio is connected, is added to the PyPSA network.
The marginal costs of this `generator' are set to the DAM prices and its nominal power is set to the grid connection capacity at the grid connection points of the portfolio.
An energy flow from the portfolio network to the `generator' represents net sales at the DAM market, a flow from the `generator' to the portfolio represents net purchases at the DAM market.
The constraint that guarantees energy conservation is automatically updated to account for this energy flow, as a standard PyPSA component is used.
The objective function $J(x)$ is also updated automatically when the `generator' is added:
\begin{equation} \label{eq:cost-da}
    J(x) \leftarrow J(x) - \sum_{t \in \mathcal{T}} c^t_\mathrm{DAM} x^t_\mathrm{DAM}\,,
\end{equation}
with $\mathcal{T}$ the time interval over which the optimization runs, $x^t_\mathrm{DAM}$ the amount of electricity sold on the DAM at time $t$ and $c^t_\mathrm{DAM}$ the DAM price at time $t$.

A bidding mechanism is modelled to restrict the interactions between portfolio and DAM.
On the German DAM, bids consist of a bid capacity and corresponding price.
If the bid is accepted, it is remunerated with the market clearing price.
This is known as a pay-as-clear or marginal pricing wholesale market.
In a MILP-formulation, either the bid capacity or the bid price can be a decision variable.
Therefore, the prices at which the bids can be made are discretized.
Bids $x^{t,j}_\mathrm{DAM,{bid}}$ can be made at price levels $c^j_\mathrm{DAM,{bid}}$ with index $j \in \mathcal{J}$.
$x^{t}_\mathrm{DAM}$ is the sum of all bids that were made by the portfolio and accepted at the market at time $t$:
\begin{equation} \label{eq:da-disp-sum}
    x^{t}_\mathrm{DAM} = \sum_{j \in \mathcal{J}} \gamma^{t,j} \cdot x^{t,j}_\mathrm{DAM,{bid}}, \quad \forall t \in \mathcal{T}\,.
\end{equation}
Here, the binary value 
\begin{equation} \label{eq:da-zuschlag}
    \gamma^{t,j} = 
    \begin{cases}
        1 & c^{j}_\mathrm{DAM,{bid}} \leq c^{t}_\mathrm{DAM} \\
        0 & \, \text{otherwise}
    \end{cases},
    \quad \forall t \in \mathcal{T}, j \in \mathcal{J}
\end{equation}
 is 1 only if the corresponding bid is accepted.

To prevent bids from cancelling each other out, all bids at a given time $t$ have to have the same sign.
This is guaranteed by adding the constraint 
\begin{equation} \label{eq:da-buy-or-sell}
    - M \cdot (1-x^{t}_\mathrm{\alpha}) \leq x^{t,j}_\mathrm{DAM,{bid}}  \leq M \cdot x^{t}_\mathrm{\alpha} \quad  \forall t \in \mathcal{T}, j \in \mathcal{J},
\end{equation}
in which $M$ is a sufficiently large number and $x^{t}_\mathrm{\alpha}$ is an additional binary decision variable, which takes value 1 only if electricity is sold on the DAM at time $t$.

\subsubsection{aFRR Market}
The aFRR bidding and market mechanisms, as well as the physical implications of the trade, that is, reserving balancing power if a bid is accepted, are modelled.
We assume that either generating units, from here on referred to as generators $g$, or BESS $b$ can contribute to the power reserves.

To model the behaviour of a generator $g$ reserving power, the following constraints are added to the model:
\begin{equation}
    x^{t,g}_\mathrm{P} + x^{t,g}_\mathrm{aFRR+} \le  P^g_\mathrm{max} \cdot x^{t,g}_\mathrm{\delta}  \quad \forall t \in \mathcal{T}, \label{eq:gen-max-power}
\end{equation}
\begin{equation}
     x^{t,g}_\mathrm{aFRR-}\le  x^{t,g}_\mathrm{P} - P^g_\mathrm{min} \cdot x^{t,g}_\mathrm{\delta} \quad \forall t \in \mathcal{T}. \label{eq:gen-min-power}
\end{equation}

Equation (\ref{eq:gen-max-power}) defines an upper bound for the positive aFRR, $x^{t,g}_\mathrm{aFRR+}$, that the generator is able to provide.
Indeed, \eqref{eq:gen-max-power} assures that the sum of the generator's power output $x^{t,g}_\mathrm{P}$ and the activated aFRR does not exceed the generator's nominal power $P_\mathrm{max}^g$.
In the same way, (\ref{eq:gen-min-power}) defines an upper bound for the negative aFRR, $x^{t,g}_\mathrm{aFRR-}$, assuring that the generator's power $x^{t,g}_\mathrm{P}$ can be decreased by up to $x^{t,g}_\mathrm{aFRR-}$, without dropping below the generator's minimal load $P_\mathrm{min}^g$.
The binary variable $x^{t,g}_\mathrm{\delta}$ is 1 only if the generator is committed at time $t$.

In order for a BESS to provide aFRR, it needs to be able to provide the full reserved power for at least one hour (according to German regulation), which is modelled by the constraints
\begin{align}
P^b_\mathrm{max} &\geq            x^{t,b}_\mathrm{aFRR+} +x^{t,b}_\mathrm{P,{out}} - x^{t,b}_\mathrm{P,{in}},&\forall t \in \mathcal{T} \label{eq:aFRR-BESS-P-pos},\\
P^b_\mathrm{max} &\geq           x^{t,b}_\mathrm{aFRR-}+x^{t,b}_\mathrm{P,{in}} - x^{t,b}_\mathrm{P,{out}},&\forall t \in \mathcal{T} \label{eq:aFRR-BESS-P-neg},\\
x^{t,b}_\mathrm{aFRR+} &\leq    x^{t-1,b}_\mathrm{E} + x^{t,b}_\mathrm{P,{in}} - x^{t,b}_\mathrm{P,{out}},&\forall t \in \mathcal{T}' \label{eq:aFRR-BESS-W-pos},\\
x^{t,b}_\mathrm{aFRR-} &\leq   E^b_\mathrm{max} - x^{t-1,b}_\mathrm{E} + x^{t,b}_\mathrm{P,{out}} - x^{t,b}_\mathrm{P,{in}} 
&\forall t \in \mathcal{T}' \label{eq:aFRR-BESS-W-neg},
\end{align}
where $\mathcal{T}'= \mathcal{T} \setminus \{t_0\}$, with $t_0$ the first element in $\mathcal{T}$. With $P^b_\mathrm{max}$ the maximum rate of (dis)charge of the BESS $b$, $x^{t,b}_\mathrm{P,{in}}$ the power flowing into the BESS and $x^{t,b}_\mathrm{P,{out}}$ the power flowing out of the BESS, \eqref{eq:aFRR-BESS-P-pos} and \eqref{eq:aFRR-BESS-P-neg} limit the (dis)charge range to $P^b_\mathrm{max}$ both with and without activation of aFRR.
Equations \eqref{eq:aFRR-BESS-W-pos} and \eqref{eq:aFRR-BESS-W-neg} make sure that the storage level $x^{t,b}_\mathrm{E}$ stays between empty and the maximum storage level $E^b_\mathrm{max}$, also in case of a one hour aFRR activation.

The aFRR market is modelled analogous to the DAM, however, aFRR markets in Germany use a pay-as-bid principle.
Therefore, the accepted bids $x^{t,k}_\mathrm{aFRR+}$ need to be weighted with their associated prices in the objective function.
In this section, only the bidding of positive aFRR is shown, negative aFRR bidding is implemented analogously.
Similar to \eqref{eq:da-disp-sum}-\eqref{eq:da-zuschlag}, the binary value 
\begin{equation*}
    \beta^{t,k}_\mathrm{aFRR+} = 
    \begin{cases}
        1 & c^{k}_\mathrm{aFRR+,{bid}} \leq c^{t}_\mathrm{aFRR+} \\
        0 &\, \text{otherwise}
    \end{cases}, \forall t \in \mathcal{T}, k \in \mathcal{K},
\end{equation*}
denotes whether a bid $x^{t,k}_\mathrm{aFRR+,{bid}}$ is accepted or rejected.
A bid is accepted if the corresponding price $c^{k}_\mathrm{aFRR+,{bid}}$ is less or equal to the marginal market price $c^{t}_\mathrm{aFRR+}$.
The accepted aFRR can be provided by generators $x^{t,g}_\mathrm{aFRR+}$ or BESSs $x^{t,b}_\mathrm{aFRR+}$, resulting in the equality constraint
\begin{equation*}
    \sum_{k \in \mathcal{K}} \beta^{t,k}_\mathrm{aFRR+} x^{t,k}_\mathrm{aFRR+,{bid}} = \sum_{g} x^{t,g}_\mathrm{aFRR+} + \sum_{b} x^{t,b}_\mathrm{aFRR+}, \forall t \in \mathcal{T}.
\end{equation*}
The aFRR bid $x^{t,k}_\mathrm{aFRR+,bid}$ with the corresponding price $c^k_\mathrm{aFRR+,{bid}}$ is added to the objective function as negative cost
\begin{multline*}
    J(x) \leftarrow J(x) 
            - \sum_{t \in \mathcal{T},k \in \mathcal{K}}\beta^{t,k}_\mathrm{aFRR+}  c^{t,k}_\mathrm{aFRR+,{bid}} x^{t,k}_\mathrm{aFRR+,bid}\\ 
            - \sum_{t \in \mathcal{T},l \in \mathcal{L}}\beta^{t,l}_\mathrm{aFRR-}  c^{t,l}_\mathrm{aFRR-,{bid}} x^{t,l}_\mathrm{aFRR-,bid}\,.
\end{multline*}

\subsection{Stochastic Optimization}
\subsubsection{Optimization Problem Formulation}
To find the optimal UC for an uncertain future, the deterministic model, as given in (\ref{eq:MILP}), is modified to:
\begin{subequations}\label{eq:stoch-MILP}
    \begin{align}
    \min_{x} &\sum_{s \in \mathcal{S}} \pi^s J(x^s), \qquad x \in \mathbb{R}^{n_\mathrm{r}\times |\mathcal{S}|} \times &\{0,1\}^{n_\mathrm{b}\times |\mathcal{S}|}\,,\\
\text{s.t. }  &A x^{s} \leq b^{s}, \qquad &\forall s \in \mathcal{S}\,,\\
    &x^{s} = 
        \begin{pmatrix}
            y \\ z^{s}
        \end{pmatrix}, \,
    x^{s'} = 
        \begin{pmatrix}
            y \\ z^{s'}
        \end{pmatrix}, \qquad &\forall s,s' \in \mathcal{S}. \label{eq:stoch-MILP-condition}
    \end{align}    
\end{subequations}
This is a scenario-based stochastic optimization, in which each scenario $s \in \mathcal{S}$ is characterized by the scenario inputs and a probability $\pi^s$.
Equation (\ref{eq:stoch-MILP}) is again a MILP.
The interpretation of (\ref{eq:stoch-MILP}) is that the objective function $J(x)$ of the original problem is replaced with the expected value of the objective function, taking into account all scenarios, and thereby the probability distributions of the optimization problem's inputs:
\begin{equation}
    \E_{s \in \mathcal{S}}[J(x^s)] = \sum_{s \in \mathcal{S}} \pi^s \, J(x^s).
\end{equation}

The constraint (\ref{eq:stoch-MILP-condition}) introduces a coupling between the different scenarios.
It states that some decision variables need to be equal for all scenarios.
For the UC problem, a single decision on the bid to make at each price should be made, under consideration of all scenarios.
To couple the decision variables for the bids, the constraints
\begin{align*} 
    x^{s,t,j}_\mathrm{DAM,{bid}} &= x^{s',t,j}_\mathrm{DAM,{bid}} \; \forall \{s,s' \in S | s \neq s'\},t \in \mathcal{T}, j \in \mathcal{J}\,
\intertext{and}
    x^{s,t,k}_\mathrm{aFRR+,{bid}} &= x^{s',t,k}_\mathrm{aFRR+,{bid}} \; \forall \{s,s' \in S | s \neq s'\},t \in \mathcal{T}, k \in \mathcal{K}
\end{align*}
are added. For negative aFRR, an analogous constraint is implemented.

To deal with the case that a bid at a certain price level is not accepted in any scenario, the following constraint is added.
The constraint forces such unaccepted bids to 0:
\begin{equation*} \label{eq:Da-unused}
    x^{s,t,j}_\mathrm{DAM,{bid}} \prod_{s' \in S} (1-\gamma^{s',t,j})= 0,  \qquad \forall s \in \mathcal{S}, t \in \mathcal{T}, j \in \mathcal{J}.
\end{equation*}

\subsubsection{Stochastic Forecasting}
When applying stochastic UC, stochastic inputs for the UC are required.
These can be stochastic price predictions for DAM or aFRR markets, but also wind, solar or district heat demand predictions, in case this is relevant for the portfolio.
If such stochastic predictions are available, scenarios can be derived from them.

As the UC tool takes scenarios as input, the method through which these scenarios are generated can be chosen freely.
When multiple variables require stochastic inputs, coupled forecasts and clustering are advisable, as this reduces the resulting number of scenarios, and hence reduces computation time.
For the example presented in Section \ref{sec:example}, simple ARIMAX prediction models were used to generate the predictions.
From the prediction-models, a big number of discrete scenarios is sampled, after which K-means clustering is used to reduce the number of scenarios over which the optimization is performed.
The probability of the K-means scenarios is given as the sum of the probabilities of the clustered scenarios.

Assuming independent variables, an ARIMAX model for a single variable can be generally written as \citep{seabold_statsmodels_2010}
\begin{equation}\label{eq:arimax}
    \begin{aligned}
        \Delta^d y^t=\phi_0  &+\phi_1 \, y^{t-1}+\phi_2 \, y^{t-2}+\ldots+\phi_p \, y^{t-p}\\
                    &+\Theta_1 \, \varepsilon^{t-1}+\Theta_2 \, \varepsilon^{t-2}+\ldots+\Theta_q \, \varepsilon^{t-q}\\
                    &+\beta \, \Psi^t+\varepsilon^t,
    \end{aligned}
\end{equation}
parameterized by $\phi$, $\Theta$ and $\beta$, under the assumption that $\varepsilon^t \sim \mathcal{N}(0,\sigma^2)$.
$\Delta$ denotes the lag operator, so that $\Delta y^t = y^t-y^{t-1}$, 
$\Delta^2 y^t=y^t-2y^{t-1}+y^{t-2}$, etc.
After integrating $d$ times, \eqref{eq:arimax} can be used to predict the time series $y$.
For prediction, the model uses $p$ latest values of $y$, $q$ latest prediction errors $\varepsilon$ and the exogenous variables $\Psi$. 

The assumption that $\varepsilon^t$ is normally distributed with a constant variance of $\sigma^2$ is justified for, e.g., district heating load.
However, prediction residues of prices often depend on the absolute value of the price, and are not normally distributed.
For variance stabilization of power market prices, $\sinh^{-1}$ can be used as a nonlinear transformation \citep{schneider_power_2011}:
\begin{equation} \label{eq:arcsinh}
    {\hat y}^{t}  = \sinh^{-1}
        \left(
            \frac{y^t-\mu}{\sigma}
        \right)
\end{equation}
$\mu$ and $\sigma$ are chosen as mean and standard deviation of $y$.
The transformation is applied before the prediction, and the inverse transformation is applied after clustering.

Typically, the time series for which forecasts are made contain seasonalities.
To forecast seasonal effects with long periods without the excessive computational effort associated with seasonal-ARIMAX models, sine and cosine functions are used to generate exogenous variables $\Psi^{T}_\mathrm{cos}$ and $\Psi^{T}_\mathrm{sin}$ for a given period $T$ \citep{smith_pmdarima_2017}:
\begin{equation}
    \Psi^{T,t}_\mathrm{cos} = \cos \left(\frac{2\pi}{T} t\right) \quad \mathrm{and} \quad
    \Psi^{T,t}_\mathrm{sin} = \sin \left(\frac{2\pi}{T} t\right)
\end{equation}

To improve prediction results, other exogenous variables can be used.
For the example in Section \ref{sec:example}, the following time series were used as exogenous variables:
\begin{itemize}
    \item the residual load forecast,
    \item the air temperature forecast in the portfolio area (used for district heating load forecasting),
    \item a weekend indicator.
\end{itemize}
The residual load forecast was calculated from the load, wind and solar power forecasts available on the ENTSO-E transparency platform. The air temperature forecast is available from the German weather service.
This selection of exogenous variables was inspired by \citet{kraft_stochastic_2022}.

To estimate the parameters of the ARIMAX model, an implementation of the Python module \textit{statsmodels} \citep{seabold_statsmodels_2010} is used.
The parameter space is searched using grid search, and optimal models are selected based on the Akaike Information Criterion (AIC).

\begin{figure}[htp!]
    \centering
    \includegraphics[width=0.7\linewidth]{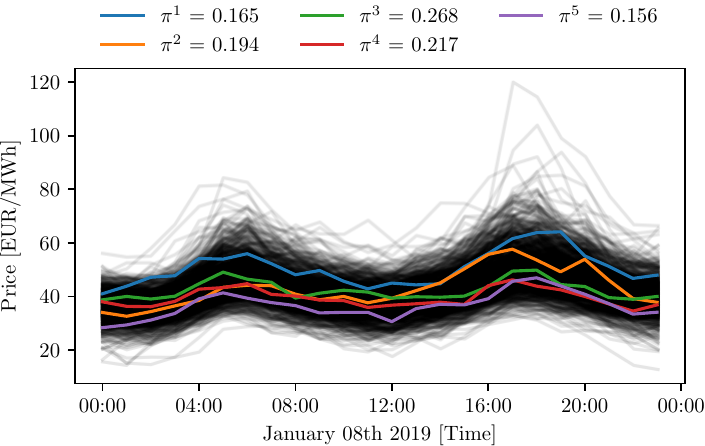}
    \caption{Clustered DAM price forecasts on 08.01.2019. The legend shows the probabilities of the different clusters. Training data: years 2019 and 2020.} \label{fig:cluster}
\end{figure}

To generate scenarios, $N_\mathrm{samp} = 1000$ realisations of the optimal ARIMAX models are generated.
K-means clustering is then used to reduce the number of scenarios to $n_\mathrm{clust} = 5$.
After the clustering, the probability $\pi^s$ of the scenario cluster $s$ is defined by
\begin{equation}
    \pi^s = \frac{N^s_\mathrm{cl}}{N_\mathrm{samp}},
\end{equation}
where $N^s_\mathrm{cl}$ is the number of scenarios represented by the scenario cluster $s$. An example of the method and its results is given for DAM price forecasts in Fig. \ref{fig:cluster}.

\subsection{Rolling Horizon and Multistaging}\label{sec:multistaging}
Using a rolling horizon, also called receding horizon, the optimal UC strategy for, e.g., the next 36h is calculated, but is only applied for the next hour or the next day, after which the optimal strategy is recalculated.
This allows portfolio managers to continuously take into account new information, e.g., new forecasts, by updating the optimization inputs.
The optimization horizon is typically larger than the time during which the strategy is applied.
This avoids final constraints or final costs from affecting the applied strategy.

\begin{figure}[htp!]
    \centering
    \includegraphics[width=0.7\linewidth,trim=4.6cm 4.35cm 4.7cm 4.37cm,clip]{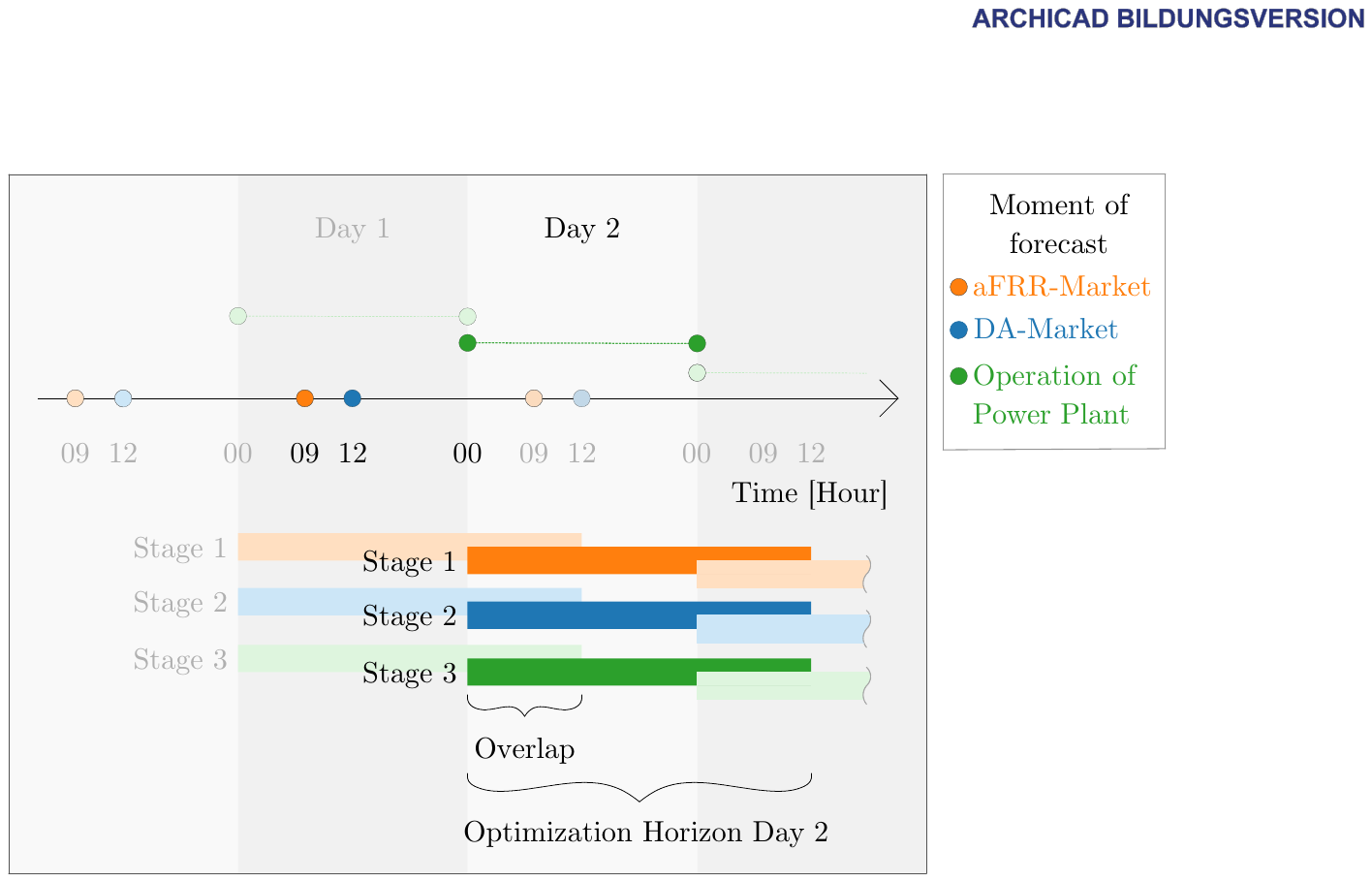}
    \caption{Overview of the stages and the rolling horizon} \label{fig:stagesOverview}
\end{figure}

A variation on rolling horizon optimization is called multistaging \citep{kraft_stochastic_2022}.
When considering the DAM and aFRR markets, bids need to be handed in by a market-closure time, after which they are accepted or rejected.
This drastically affects the uncertainty, as well as the amount of decision variables, in the optimization problem.
As a result, some `stages' of the rolling horizon are less complex optimization problems than others.
The market-closure times and the associated stages are shown schematically in Fig.~\ref{fig:stagesOverview}.

The following stages are defined:
\begin{itemize}
    \item \textbf{Stage 1} is run just before the aFRR market-closure time, which is at 9 a.m. on the last day preceding the allocation day (for the German spot market).
    This stage has the largest temporal distance to the allocation day, hence the highest uncertainty. 
    The aFRR prices, DAM prices, as well as wind, solar, and heating demand (if relevant for the portfolio) are uncertain.
    \item \textbf{Stage 2} of the optimization is run just before DA market-closure, which is at 12 p.m. for the German spot market.
    The aFRR bids have been made, and the auction results are available.
    Therefore, the aFRR price forecast is replaced with the available auction prices, and the aFRR bids are fixed to their values obtained in Stage 1.
    Updated forecasts for remaining uncertain quantities, most notably the DAM prices, are used.
    \item \textbf{Stage 3} is run at the beginning of the allocation/delivery day. 
    At this time, all auction results are known.
    In this stage, the UC optimization only needs to find the optimal use of the available assets to fulfill the contracts.
    If large uncertainty would still be present at the beginning of Stage 3, e.g., in the availability of wind or solar assets, one would run this stage in a rolling horizon fashion, regularly updating forecasts and UC.
    In our implementation, a single deterministic optimization is run for Stage 3, assuming that all uncertainty is resolved at this time.
    This can, however, easily be changed.
\end{itemize}

\section{Application Example} \label{sec:example}
The methods described in Section \ref{sec:methods} can be used individually or combined to perform UC.
As proof-of-concept, a small portfolio consisting of a waste-to-energy (WtE) plant, a BESS and a thermal heat storage system is optimized to trade on the DAM and aFRR markets and continuously cover a required district heating demand.
Specifications from an existing $\qty{60}{\mathrm{MW}_\mathrm{th}}$ WtE plant and connected district heating network in Germany are used.
The plant can deliver up to \qty{30}{MW} of heat and produce up to \qty{15}{MW} of electricity.
It should be noted that, unlike most thermal plants that incur fuel costs, WtE facilities earn income by disposing of waste.
In normal operation, the plant is the only heat provider for the district heating network.
A notional \qty{6}{MWh}/\qty{6}{MW} BESS and \qty{50}{MWh}/\qty{10}{MW} heat storage system are added to allow flexible operation of the portfolio.
DA, aFRR and district heat demand data of a two-year period are used to estimate the parameters $\phi$, $\Theta$ and $\beta$ of the ARIMAX forecasting models \eqref{eq:arimax} which can be found in \citet{welfonder_pypsa-stochuc_2024}.
For simplicity’s sake, independent variables are assumed.
For each optimization stage, five scenarios are generated for each uncertain input, and are combined combinatorially, following the method outlined above.

\subsection{Applying Power and Energy Markets}

\begin{figure}
    \centering
    \includegraphics[width=0.7\linewidth]{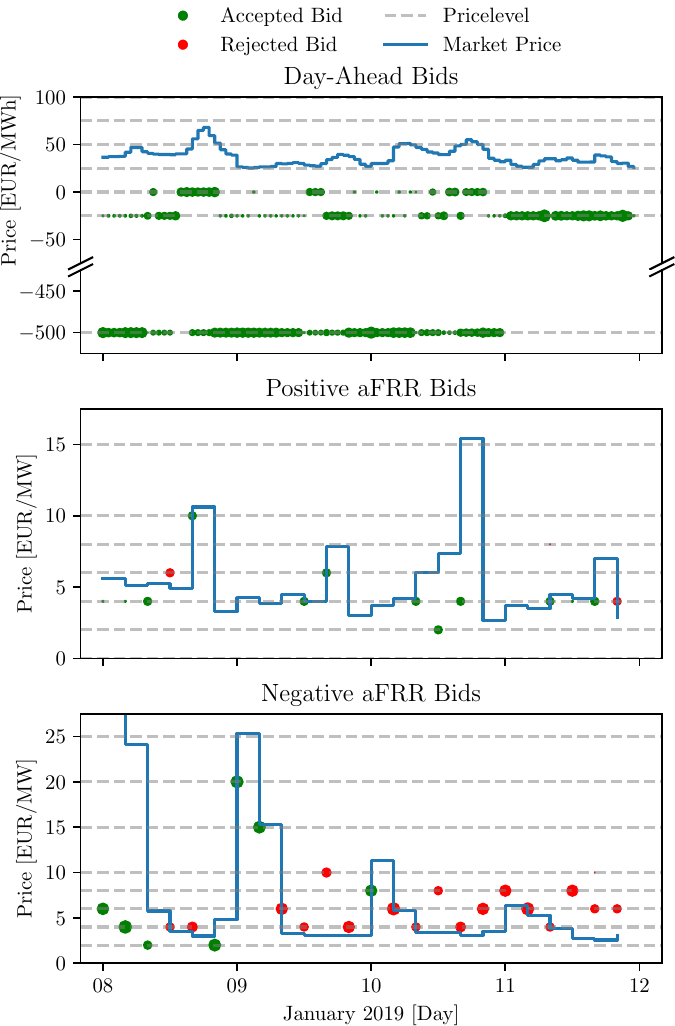}
    \caption{Bids submitted to the DAM (after Stage 2) and aFRR markets (after Stage 1) and the market prices. The size of the dots reflects the size of the bid (in MWh), the colour indicates whether the bid was accepted or not.} \label{fig:bids}
\end{figure}

\begin{figure*}
    \centering
    \includegraphics[width=\textwidth]{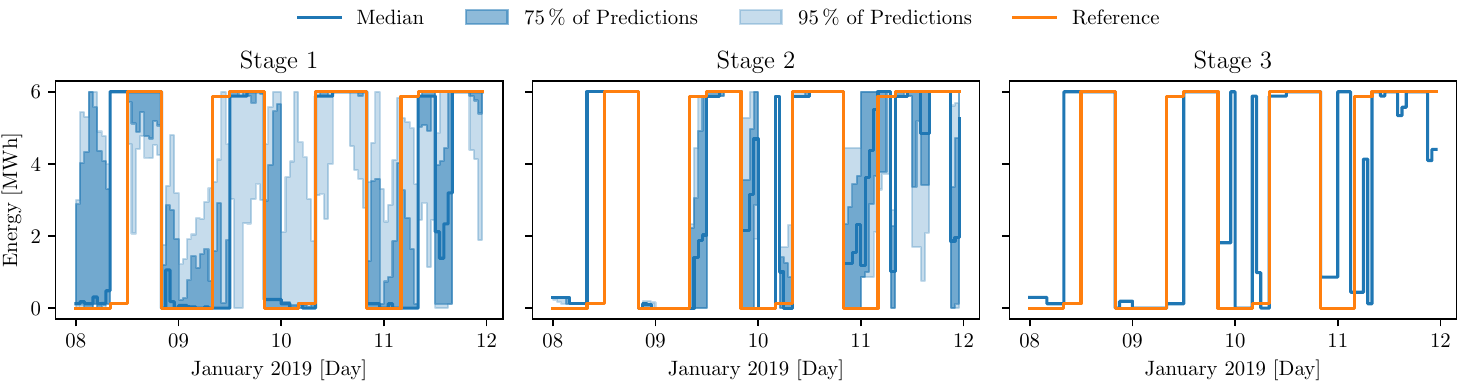}
    \caption{Stochastic UC optimization results of different stages for BESS storage level. As reference, the result of a deterministic UC optimization is shown.} \label{fig:BESS-soc}
\end{figure*}

When implementing DA and aFRR markets, one of the outputs of the UC are bids that are submitted to these markets, as is seen in Fig.~\ref{fig:bids}.
The different market mechanisms (pay-as-clear for DAM and pay-as-bid for aFRR) result in different bidding strategies seen in the figure. On the aFRR market, bids are placed around the expected market prices to maximize revenues.
This strategy strongly relies on the forecasting accuracy, as wrong forecasts might result in too high bids that are rejected (as is seen, e.g., for many negative aFRR bids).
On the DAM, bids are placed around the plant's (negative) marginal costs, and the minimal bid price (\qty{-500}{EUR/MWh}).
This reflects the pay-as-clear mechanism of the market, where revenue does not increase with higher bids.
Instead, bids would typically be placed at the marginal cost level.
The bids at the minimal bid price are a consequence of accepted aFRR bids: when aFRR is allocated, the portfolio must have committed assets, often requiring active power delivery to the grid.

\subsection{Applying Stochastic Optimization}

Fig.~\ref{fig:BESS-soc} shows the UC results for the BESS storage level, grouped per optimization stage. 
The results of a deterministic UC with perfect forecast are shown as reference.
For the stochastic optimization, the median and areas encompassing \num{75} and \qty{95}{\%} of the predictions are shown.
The effect of the stochastic optimization can be seen best in the first stage, where uncertainty is highest, as it is not yet known which DAM and aFRR bids will be accepted.
This translates in a large uncertainty of the optimal storage levels, which depend on the forecast scenario.

\subsection{Applying Rolling Horizon and Multistaging}
In Fig. \ref{fig:BESS-soc}, multiple optimization stages are used.
It is seen that, as the aFRR market closes between Stage 1 and Stage 2 and the DA market closes between Stage 2 and Stage 3, uncertainty systematically becomes less, and the results show more certainty on the optimal UC and BESS operation. 
When comparing Fig. \ref{fig:bids} and Stage 3 in Fig. \ref{fig:BESS-soc}, it is seen that the BESS discharges at DAM price peaks, and is often fully charged at times when positive aFRR is offered, and empty when negative aFRR is offered.

The sharp changes in optimal state of charge that are seen at some transitions between days in Stage 2 and Stage 3 of Fig.~\ref{fig:BESS-soc} are caused by the rolling horizon approach, for which new optimizations, with updated forecasts, are executed at midnight.
Because bids are submitted per day, the effect of the new forecasts and updated optimal UC is visible at the change between days.

\section{Outlook}
The open source stochastic unit commitment tool presented in this paper can be used flexibly to solve the UC problem for various portfolios.
Currently, aFRR and DA markets are implemented according to the German market mechanisms.
We invite the scientific community to implement market structures of other countries, or to combine the tool with other forecasting methods, which could result in better bidding on the negative aFRR market.
In order to leverage the compatibility with PyPSA, the tool should be kept in line with the latest PyPSA standards.

\section{Acknowledgements}
This research is part of the project ``CHP 4.0 -- Regional combined heat and power plants in a changing energy system.'' It is funded by the German Federal Ministry for Economic Affairs and Climate Action under grant number FKZ 03EE5031 E.

\bibliographystyle{unsrtnat} 
\bibliography{bibliography}
\end{document}

%% file: abstract.tex
This paper presents an open source stochastic unit commitment (UC) optimization tool, which is available on GitHub. In addition, it presents an example use case in which UC optimization is done for a waste-to-energy plant with heat storage and a battery energy storage system (BESS) in Germany, under uncertain day-ahead and balancing power (aFRR) market prices as well as heat load uncertainty. The tool consists of multiple modular extensions for the Python for Power System Analysis (PyPSA) framework, namely the implementation of market and bidding mechanisms, stochastic optimization and multistaging.